\documentclass[runningheads]{svjour3}
\smartqed  % flush right qed marks, e.g. at end of proof
\usepackage{amssymb,euscript}
\usepackage{graphicx}
\journalname{Science China Mathematics}

\def\dfrac#1#2{{\displaystyle {{#1}}\over{\displaystyle{#2}}}}
\begin{document}

\title{Error bounds of a quadrature formula with multiple nodes for the Fourier-Chebyshev coefficients for analytic functions}
 \thanks{
 %Grants or other notes
}

\titlerunning{Error bounds of a quadrature formula for the Fourier-Chebyshev coefficients}        % if too long for running head

\author{Aleksandar V. Pej\v cev \and Miodrag M. Spalevi\' c % \and %etc.
}


\institute{ %\at
           Department of Mathematics, University of Beograd, \\
              Faculty of Mechanical Engineering, Kraljice Marije 16,\\
              11120 Belgrade 35, Serbia\\
              \email{apejcev@mas.bg.ac.rs, mspalevic@mas.bg.ac.rs}           %
              \\
           %\and
           %Miodrag M. Spalevi\' c \at
           %Department of Mathematics, University of Beograd, \\
            %  Faculty of Mechanical Engineering, Kraljice Marije 16,\\
             % 11120 Belgrade 35, Serbia\\
             % \email{mspalevic@mas.bg.ac.rs}           %
              %\\
}

\date{Received: 31 October, 2016 / Revised: 30 December, 2017}

\maketitle

\begin{abstract} Three kinds of effective error bounds of the quadrature formulas with multiple nodes that are  generalizations of the well known Micchelli-Rivlin quadrature
formula, when the integrand is a function analytic in the regions bounded by confocal ellipses, are given. A numerical example which illustrates the calculation of these error bounds is included.

\vskip 0.1in

\noindent {\bf [This paper has been accepted for publication in SCIENCE CHINA Mathematics.]}

\end{abstract}

%Insert your abstract here. Include up to five keywords.

\section{Introduction}

We consider the quadrature formula with multiple nodes
\begin{equation}\label{glavnaqf}
\int_{-1}^1 f(t)T_n(t)\,\frac{dt}{\sqrt{1-t^2}}=\sum_{\nu=1}^n\sum_{i=0}^{2s-1}A_{i,\nu}f^{(i)}(\xi_\nu)+R_{n,s}(f)
\end{equation}
for calculating the Fourier-Chebyshev coefficients of an analytic function $f$ ($n\in\mathbb{N}, s\in\mathbb{N}$), with respect to the Chebyshev weight function of the first kind $\omega(t)=1/\sqrt{1-t^2}$. $T_n$ is the Chebyshev polynomial of the first kind of degree $n$,
\[
T_n(t)=\cos(n\arccos\,
t)=2^{n-1}\,(t-\xi_1)\cdots(t-\xi_n),\quad t\in(-1,1).
\]

 The quadrature formula (\ref{glavnaqf}) has been firstly mentioned in \cite[p. 383]{bojanov}, and then analyzed in more details in \cite{GVM_RO_MMS}. It has the algebraic degree of precision $n(2s+1)-1$. Its special case $s=1$ represents the well-known Micchelli-Rivlin quadrature formula introduced in \cite{mr}. Micchelli and Rivlin \cite{mr} considered a quadrature formula of the highest algebraic degree of precision for the Fourier-Chebyshev coefficients $a_n(f)$,
\[
a_n(f)=\int_a^b T_n(t)f(t)\omega(t)\,dt,
\]
which is based on the divided differences
of $f'$ at the zeros of the Chebyshev polynomial $T_n$.
For more details on this subject see \cite{boj},
\cite{bojanov}, \cite{rdevore}, \cite{mr10}, \cite{MSmc2013}.

In \cite{APejcevMMS} we considered the error bounds of the Micchelli-Rivlin quadrature formula for analytic functions. In this paper we will consider the corresponding error bounds of its generalizations (\ref{glavnaqf}) ($s>1$).

\section{Error bounds of the quadrature formula (\ref{glavnaqf}) for analytic functions}

Let $\Gamma$ be a simple closed curve in the complex plane surrounding $[-1,1]$ and $\mathcal{D}$ its interior.
Let $f$ be an analytic function in $\mathcal{D}$ and
continuous on $\overline{\mathcal{D}}$.
 If the values
of the function $f$ and of its derivatives $f^{(i)},\ i=1,\dots,2s-1$  in the
nodes $x_1,x_2,...,x_n(\in[-1,1])$ are known, then the error
of Hermite interpolation of the function $f$ can be written in the
form (see Gon\v carov \cite{gon})
\begin{equation}\label{kvadratura}
r_{n,s}(f;t)=f(t)-\sum_{\nu=1}^n\sum_{i=0}^{2s-1}
\ell_{i,\nu}(t)f^{(i)}(x_{\nu})=\dfrac{1}{2\pi
i}\oint_{\Gamma}\dfrac{f(z)\Omega_{n,s}(t)}{(z-t)\Omega_{n,s}(z)}\,dz,
\end{equation}
where $\ell_{i,\nu}$ are the fundamental polynomials of the Hermite
interpolation and $\Omega_{n,s}(z)=\prod_{\nu=1}^n(z-x_{\nu})^{2s}$.

If we choose $x_{\nu}$ to be the zeros of the Chebyshev polynomial
of the first kind, i.\,e., $x_{\nu}=\xi_{\nu}$, after multiplying by
(\ref{kvadratura}) with $\omega(t)T_n(t)$, where
$\omega(t)={1}/{\sqrt{1-t^2}}$, and integrating in $t$ over
$(-1,1)$, we get a contour integral representation of the
remainder term in (\ref{glavnaqf}).

We get the representation
\begin{equation}\label{ostatak}
R_{n,s}(f)\equiv R_{n,s}(fT_n)=\frac{1}{2\pi i}\oint_{\Gamma}K_{n,s}(z)f(z)dz,
\end{equation}
where the {\it kernel} is given by
\begin{equation}\label{jezgro}
K_{n,s}(z)=\dfrac{\rho_{n,s}(z)}{T_n^{2s}(z)},
\end{equation}
and
\begin{equation}\label{ron}
\rho_{n,s}(z)=\int_{-1}^1\dfrac{\omega(t)}{z-t}\,T_n^{2s+1}(t)dt.
\end{equation}

From (\ref{ostatak}) we obtain the error bound
\begin{equation}\label{ocena}
\left|R_{n,s}(fT_n)\right|\le
\frac{\ell(\Gamma)}{2\pi}\left(\max_{z\in\Gamma}|K_{n,s}(z)|\right)
\left(\max_{z\in\Gamma}|f(z)|\right),
\end{equation}
where $\ell(\Gamma)$ is the length of the contour $\Gamma$.

More in general, if we apply the H\"{o}lder inequality to (\ref{ostatak}), we get
\begin{eqnarray*}
\left|R_{n,s}(fT_n)\right|&=&\dfrac{1}{2\pi }\left|\oint_{\Gamma}K_{n,s}(z)f(z)dz\right|\\
&\leq& \dfrac{1}{2\pi }\left(\oint_{\Gamma}|K_{n,s}(z)|^r|dz|\right)^{1/r}\left(\oint_{\Gamma}|f(z)|^{r'}|dz|\right)^{1/{r'}},\\
\end{eqnarray*}
i.e.
\begin{equation}\label{Helder}
\|R_{n,s}(fT_n)\|\leq\dfrac{1}{2\pi } \|K_{n,s}\|_r\|f\|_{r'},
\end{equation}
where $1\leq r \leq +\infty$, $1/r+1/r'=1$ and
\[
\|f\|_r=\left\{\begin{array}{llr}
\left(\displaystyle\oint_{\Gamma}|f(z)|^{r}|dz|\right)^{1/{r}}, &1\leq r < +\infty, \\
\displaystyle\max_{z\in\Gamma}|f(z)|, &r=+\infty.
\end{array}\right.
\]

In the case $r=+\infty$, $r'=1$, the estimate (\ref{Helder}) reduces to
\begin{equation}
|R_{n,s}(fT_n)|\leq
\dfrac{1}{2\pi}\left(\max_{z\in\Gamma}|K_{n,s}(z)|\right)\left(\oint_{\Gamma}|f(z)||dz|\right),
\end{equation}
which leads to the error bound (\ref{ocena}) (see, e.\,g., \cite{gauvar},
\cite{sch}, \cite{geno}, \cite{NM2016PejSpa}). We refer to it as the
$L^\infty$-error bound.

 On the other side, for $r=1$ ($r'=+\infty$) the estimate
(\ref{Helder}) reduces to
\begin{equation}\label{eocena}
\left|R_{n,s}(fT_n)\right|\le
\frac{1}{2\pi}\left(\oint_{\Gamma}|K_{n,s}(z)||dz|\right)
\left( \max_{z\in\Gamma}|f(z)|\right),
\end{equation}
which is evidently stronger than (\ref{ocena}) because of the inequality
\begin{equation}
\oint_{\Gamma}\left|K_{n,s}(z)\right|\,|dz|\leq
\ell(\Gamma)\left(\max_{z\in \Gamma}|K_{n,s}(z)|\right).
\end{equation}

We refer to (\ref{eocena})  as to the $L^1$-error bound.

 In this paper we take $\Gamma=\mathcal{E}_\rho$, where the ellipse
$\mathcal{E}_\rho$ is given by

\begin{equation}\label{elipsa}
\mathcal{E}_\rho=\left\{z\in \mathbb{C}\, \left|\right.\ z=\frac
12\left(u+u^{-1}\right),\ 0\le\theta\le
2\pi\right\},\quad u=\rho\,e^{i\theta}.
\end{equation}

The choice of the family of ellipses $\mathcal{E}_\rho$ as basic contours of integration  is natural when dealing with analytic functions in a neighborhood of $[-1,1]$, since they are the level curves of the Green Function of $\mathbb{C}\setminus[-1,1]$ with pole at infinity, in such a way that for $\rho\to1^+$,  $\mathcal{E}_\rho$ tends to $[-1,1]$ and $\rho\to\infty$, the interior of $\mathcal{E}_\rho$ approaches the whole complex plane (interesting when dealing with entire integrands, as in Section \ref{sect6}).

\section{$L^\infty$-error bounds based on the analysis of the maximum modulus of the kernel}

 We have from (\ref{ron}), by substituting $t=\cos\theta$,
\begin{eqnarray*}
\rho_{n,s}(z)&=&\int_0^{\pi}\dfrac{[\cos{n\theta}]^{2s+1}}{z-\cos{\theta}}d\theta\\
&=&\dfrac{1}{2^{2s}}\int_0^{\pi}\dfrac{1}{z-\cos{\theta}}\left(\sum_{k=0}^{s}{2s+1 \choose k}\cos{(2s+1-2k)n\theta}\right)d\theta,
\end{eqnarray*}
where we used \cite[Eq. 1320.7]{Ryzhik}. Now the kernel has the form
\[
K_{n,s}(z)=\dfrac{\dfrac{1}{2^{2s}}\sum_{k=0}^{s}{2s+1 \choose k}\int_0^{\pi}\dfrac{\cos{(2s+1-2k)n\theta}}{z-\cos{\theta}}d\theta}{[T_n(z)]^{2s}},
\]
i.\,e.
\[
K_{n,s}(z)=\dfrac{\dfrac{1}{2^{2s}}\sum_{k=0}^{s}{2s+1 \choose k}\dfrac{\pi}{\sqrt{z^2-1}}\left(z-\sqrt{z^2-1}\right)^{2s+1-k}}{[T_n(z)]^{2s}},
\]
where we used (see, e.\,g., \cite{gauvar})
\[
\int_0^{\pi}\dfrac{\cos{m\theta}}{z-\cos{\theta}}d\theta=\dfrac{\pi}{\sqrt{z^2-1}}\left(z-\sqrt{z^2-1}\right)^m, \ m\in \mathbb{N}_0.
\]
Substituting $z=\frac{1}{2}(u+u^{-1})$ $\left(u=z+\sqrt{z^2-1}\right)$, using
\begin{equation}\label{cebpol}
T_n(z)=\left(u^n+u^{-n}\right)/{2},
\end{equation}
we get
\begin{eqnarray*}
K_{n,s}(z)&=&\dfrac{\dfrac{1}{2^{2s}}\sum_{k=0}^{s}{2s+1 \choose k}\dfrac{2\pi}{\left(u-u^{-1}\right)}u^{k-2s-1}}{\left[\left(u^n+u^{-n}\right)/{2}\right]^{2s}}\\
&=&\dfrac{2\pi\sum_{k=0}^{s}{2s+1 \choose k}\left(u^{2n}\right)^k}{u^{(2s+1)n}\left(u^n+u^{-n}\right)^{2s}\left(u-u^{-1}\right)}.
\end{eqnarray*}
With the usual notation (see \cite{gauvar})
\begin{equation}\label{aiovi}
a_j=a_j(\rho)=\dfrac{1}{2}(\rho^j+\rho^{-j}), \ j\in
\mathbb{N}\quad (\rho>1),
\end{equation}
when $u=\rho e^{i\theta}$, we have
\begin{eqnarray*}
a&=&\left|\sum_{k=0}^{s}{2s+1 \choose k}\left(u^{2n}\right)^k\right|^2\\
&=&\left(\sum_{k=0}^{s}{2s+1 \choose k}\rho^{2nk}\cos{2nk\theta}\right)^2+\left(\sum_{k=0}^{s}{2s+1 \choose k}\rho^{2nk}\sin{2nk\theta}\right)^2,\\
\left|u-u^{-1}\right|^2&=&2(a_2-\cos{2\theta})=2b,\\
\left|u^n+u^{-n}\right|^2&=&2(a_{2n}+\cos{2n\theta})=2c,
\end{eqnarray*}
and
\begin{equation}\label{kvadratmodula}
|K_{n,s}(z)|^2=\dfrac{\pi^2}{2^{2s-1}\rho^{2(2s+1)n}}\cdot\dfrac{a}{bc^{2s}}.
\end{equation}

Let us denote by $A,B,C$ the values of $a,b,c$ at $\theta=0$, respectively.

Now we can formulate the main statement.
\begin{theorem}\label{teorema1} For each fixed $n\in\mathbb{N}$ there exists $\rho_0=\rho_0(n)$ such that
\[
\max_{z\in \mathcal{E}_{\rho}}\left|K_{n,s}(z)\right|=
\left|K_{n,s}\left(\frac{1}{2}(\rho+\rho^{-1})\right)\right|,
\]
for each $\rho>\rho_0$.
\end{theorem}

\noindent {\bf Proof.} This condition is equivalent to
\[
\dfrac{a}{bc^{2s}}\leq\dfrac{A}{BC^{2s}},
\]
i.\,e.
\[
I=aBC^{2s}-Abc^{2s}\leq 0,
\]
for each $\rho$ greater than some $\rho_0$ on the domain $(1,+\infty)$. The member with the highest degree of $\rho$ in this expression is
\begin{eqnarray*}
&&{2s+1 \choose s}^2\rho^{4ns}\cdot(-1)\cdot\left(\dfrac{1}{2}\rho^{2n}\right)^{2s}-{2s+1 \choose s}^2\rho^{4ns}\cdot(-\cos{2\theta})\cdot\left(\dfrac{1}{2}\rho^{2n}\right)^{2s}\\
&=&\dfrac{1}{2^{2s}}\left(\cos{2\theta}-1\right){2s+1 \choose s}^2\rho^{8ns},
\end{eqnarray*}
and it is obviously negative for each $\theta\in(0,\pi]$. \hfill\hfill \qed

The empirical results show that we can take $\rho_0=1$ in almost all the cases.

\section{Error bounds based on an expansion of the remainder term}

If $f$ is an analytic function in the interior of
$\mathcal{E}_{\rho}$, it has the expansion
\begin{equation}\label{cebred}
f(z)=\sum_{k=0}^{\infty}{'}\alpha_k T_k(z),
\end{equation}
where $\alpha_k$ are given by
\[
\alpha_k=\dfrac{1}{\pi}\int_{-1}^1(1-t^2)^{-1/2}f(t)T_k(t)dt.
\]
The series (\ref{cebred}) converges for each $z$ in the interior
of $\mathcal{E}_{\rho}$. The prime in the corresponding sum denotes
that the first term is taken with the factor ${1}/{2}$.

\begin{lemma}\label{Lema1}
If $z\notin [-1,1]$, then, the following expansion holds
\begin{equation}\label{razvojceb}
\dfrac{1}{[T_n(z)]^{2s}}=\sum_{k=0}^{+\infty}\beta^{(s)}_{n,k} u^{-2ns-k},
\end{equation}
where
\begin{equation}\label{bete}
\displaystyle\beta^{(s)}_{n,k}=\left\{\begin{array}{ll}
\displaystyle 2^{2s}(-1)^j {{j+2s-1} \choose {2s-1}}, & k=2jn, \\
 0, & \mbox{otherwise}.
\end{array}\right.
\end{equation}
\end{lemma}

\noindent{\bf Proof.}
We know that if $x\in \mathbb{C}$, $|x|<1$, then
\begin{equation}
\dfrac{1}{(1-x)^{\nu+1}}=\sum_{k=\nu}^{+\infty}{k \choose \nu}x^{k-\nu}=\sum_{j=0}^{+\infty}{j+\nu \choose \nu}x^{j} \qquad (\nu=0,1,2,...).
\end{equation}

Using this fact and (\ref{cebpol}), with $u=\rho e^{i\theta}$,
$\rho>1$, $z= (u+u^{-1})/{2}$, we get
\begin{eqnarray*}
\dfrac{1}{[T_n(z)]^{2s}}&=&\left[\dfrac{1}{2}(u^n+u^{-n})
\right]^{-2s}=2^{2s}u^{-2ns}\left(\dfrac{1}{1-\left(-u^{-2n}\right)}\right)^{(2s-1)+1}\\
&=& 2^{2s}\sum_{j=0}^{+\infty}(-1)^j{{j+2s-1} \choose {2s-1}}u^{-2ns-2nj},
\end{eqnarray*}
which completes the proof.\hfill\hfill\qed

\begin{lemma}\label{Lema2}
If $z\notin [-1,1]$, $\rho_{n,s}$ can be expanded as
\begin{equation}\label{razvojrho}
\rho_{n,s}(z)=\sum_{k=0}^{+\infty}\gamma_{n,k}^{(s)}u^{-n-k-1},
\end{equation}
where
\begin{equation}\label{gamma}
{\gamma}_{n,k}^{(s)}=\left\{\begin{array}{ll}
 \displaystyle\frac{\pi}{2^{2s-1}}\displaystyle\sum_{\nu=0}^j{{2s+1} \choose{s-\nu}}, & k=2nj,2nj+2,\dots,2n(j+1)-2, \  j\in\mathbb{N}_0, \\[0.1in]
 0, & \mbox{\rm otherwise}. \\
\end{array}\right.
\end{equation}
\end{lemma}

{\bf Proof.}
It is obvious that we have the same situation with those coefficients as in \cite{Baza} and the statement directly follows from \cite{filomat}.\qed

Now, substituting (\ref{razvojceb}) and (\ref{gamma}) in
(\ref{jezgro}), we obtain
\begin{equation}\label{razvojjezgro}
K_{n,s}(z)=\sum_{k=0}^{+\infty}\omega_{n,k}^{(s)}u^{-(2s+1)n-k-1},
\end{equation}
where
\begin{equation}\label{konvolucija}
\omega_{n,k}^{(s)}=\sum_{j=0}^k\beta_{n,j}^{(s)}\gamma_{n,k-j}^{(s)}.
\end{equation}

\begin{theorem}
The remainder term $R_{n,s}(f)$ can be represented in the form
\begin{equation}\label{razvojostatak}
R_{n,s}(f)=\sum_{k=0}^{+\infty}\alpha_{(2s+1)n+k}\ \epsilon_{n,k}^{(s)},
\end{equation}
where the coefficients $\epsilon_{n,k}^{(s)}$ are independent on $f$.
Furthermore, if $f$ is an even function then $\epsilon_{n,2j+1}=0$
$(j=0,1,...).$
\end{theorem}

\noindent{\bf Proof.}
By substituting (\ref{cebred}) and (\ref{razvojjezgro}) in
(\ref{ostatak}) we obtain
\begin{eqnarray*}
R_{n,s}(f)&=&\dfrac{1}{2\pi i}\int_{\mathcal{E}_\rho}\left(\sum_{k=0}^{\infty}{'}\alpha_k T_k(z)\sum_{k=0}^{+\infty}\omega_{n,k}^{(s)}u^{-(2s+1)n-k-1}\right)dz \\
&=&\sum_{k=0}^{+\infty}\left(\dfrac{1}{2\pi i}\sum_{j=0}^{+\infty}{'}\alpha_j \right.
\left.\int_{\mathcal{E}_\rho}T_j(z)u^{-(2s+1)n-k-1}dz\right)\omega_{n,k}^{(s)}.
\end{eqnarray*}
Applying Lemma 5 from \cite{Hunter}, this reduces to
(\ref{razvojostatak}) with
\begin{equation}\label{veza}
\epsilon_{n,0}^{(s)}=\dfrac{1}{4}\omega_{n,0}^{}(s), \
\epsilon_{n,1}^{(s)}=\dfrac{1}{4}\omega_{n,1}^{(s)}, \
\epsilon_{n,k}^{(s)}=\dfrac{1}{4}(\omega_{n,k}^{(s)}-\omega_{n,k-2}^{(s)}), \
k=2,3,...\ .
\end{equation}
When $k$ is odd, since $\omega(t)=\omega(-t)$ it follows from
(\ref{konvolucija}) and Lemmas \ref{Lema1} and \ref{Lema2} that $\omega_{n,k}^{(s)}=0$,
and hence $\epsilon_{n,k}^{(s)}=0$.\hfill\hfill\qed

\subsection{Error bounds based on the estimation of the coefficients}
In general, the Chebyshev-Fourier coefficients $\alpha_k$ in
(\ref{cebred}) are unknown.
 However, Elliot \cite{Eliot} described a number of ways of estimating or bounding them. In particular, under our assumptions
\begin{equation}\label{cebocena}
|\alpha_k|\leq \dfrac{2}{\rho^k}\left(\max_{z\in \mathcal{E}_{\rho}}|f(z)| \right).
\end{equation}
By using (\ref{bete}), (\ref{gamma}), (\ref{konvolucija}), if and
only if $k=2jn$, $j\in \mathbb{N}_0$, we have
\begin{eqnarray*}
\omega_{n,2jn}^{(s)}&=&\beta_0^{(s)}\gamma_{2jn}^{(s)}+\beta_{2n}^{(s)}\gamma_{(2j-2)n}^{(s)}+...+\beta_{(2j-2)n}^{(s)}\gamma_{2n}^{(s)}+\beta_{2jn}^{(s)}\gamma_0^{(s)},\\
\omega_{n,2jn-2}^{(s)}&=&\beta_0^{(s)}\gamma_{2jn-2}^{(s)}+\beta_{2n}^{(s)}\gamma_{(2j-4)n}^{(s)}+...+\beta_{(2j-2)n}^{(s)}\gamma_{2n-2}^{(s)},
\end{eqnarray*}
which implies (cf. \cite{filomat})
\begin{equation}\label{epsilon}
\epsilon_{n,k}^{(s)}=\left\{\begin{array}{ll}
\displaystyle \pi\sum_{j=m-s}^m(-1)^j{{j+2s-1}\choose{2s-1}}{{2s+1}\choose {s-(m-j)}}, & k=2nm, \ m\in\mathbb{N}_0, \\[0.in]
 0, & \mbox{\rm{otherwise}}. \\
\end{array}\right.
\end{equation}
The last sum can be rewritten in the form
\begin{equation}\label{epsilon1}
(-1)^{m-s}\pi\sum_{i=0}^{s}(-1)^{i}{m+s-1+i \choose 2s-1}{2s+1 \choose i}.
\end{equation}

Now we can formulate and prove the following statement.

\begin{lemma}
For each $t\in\mathbb{N}_0$, it holds
\begin{equation}\label{pom}
\sum_{i=0}^{t}(-1)^{i}{m+s-1+i\choose 2s-1}{2s+1\choose i}=(-1)^{t}\dfrac{s(2m+2s+2)-t}{(m+s)(m+s+1)}{m+s+t \choose 2s}{2s \choose
t}.
\end{equation}
\end{lemma}

\noindent{\bf Proof}.
We will prove this using the mathematical induction principle over $t$.
For $t=0$ we need to prove
\[
{m+s-1 \choose 2s-1}{2s+1 \choose 0}=\dfrac{s(2m+2s+2)}{(m+s)(m+s+1)}{m+s \choose
2s}{2s \choose 0},
\]
which is obvious.

If we suppose that (\ref{pom}) holds for some $t\in
\mathbb{N}_0$ and we want to deduce that it holds for $t+1$, we have to confirm the identity
\begin{eqnarray*}
&&(-1)^{t}\dfrac{s(2m+2s+2)-t}{(m+s)(m+s+1)}{m+s+t \choose 2s}{2s \choose
t}+{(-1)}^{t+1}{m+s+t\choose 2s-1}{2s+1 \choose t+1}\\
&&=(-1)^{t+1}
\dfrac{s(2m+2s+2)-t-1}{(m+s)(m+s+1)}{m+s+t +1\choose 2s}{2s \choose
t+1},
\end{eqnarray*}
i.\,e.
\begin{eqnarray*}
&&\dfrac{s(2m+2s+2)-t}{(m+s)(m+s+1)}{m+s+t \choose 2s}{2s \choose
t}\\
&&+\dfrac{s(2m+2s+2)-t-1}{(m+s)(m+s+1)}{m+s+t +1\choose 2s}{2s \choose
t+1}\\
&&={m+s+t\choose 2s-1}{2s+1 \choose t+1},
\end{eqnarray*}
i.\,e.
\begin{eqnarray*}
&&\dfrac{s(2m+2s+2)-t}{(m+s)(m+s+1)}\cdot\dfrac{m-s+t+1}{2s}{m+s+t \choose 2s-1}\cdot\dfrac{t+1}{2s+1}{2s+1 \choose
t+1}\\
&&+\dfrac{s(2m+2s+2)-t-1}{(m+s)(m+s+1)}\cdot\dfrac{m+s+t+1}{2s}{m+s+t \choose 2s-1}\cdot\dfrac{2s-t}{2s+1}{2s+1 \choose
t+1}\\
&&={m+s+t\choose 2s-1}{2s+1 \choose t+1},
\end{eqnarray*}
which is equivalent to
\begin{eqnarray*}
&&\left(s(2m+2s+2)-t \right)\left(m-s+t+1\right)(t+1)\\
&&+\left(s(2m+2s+2)-t-1\right)\left(m+s+t+1\right)(2s-t)\\
&&=2s(2s+1)(m+s)(m+s+1).
\end{eqnarray*}
We can directly confirm this identity, but we can also do it in a little bit shorter way. Namely, the left-hand side presents the polynomial in $t$ and its degree is less than 3 (there, no power of $t$ higher than $3$ appears, and the corresponding coefficient is equal to $(-1)\cdot 1 \cdot 1+(-1)\cdot 1\cdot (-1)=0$ ), and then it is enough to show that the identity holds for three different values of $t$. The easiest choice of those values would be $t_1=-1$, $t_2=2s$ and $t_3=0$ and in each of them we only have to show the equality of two products.
\hfill\hfill\qed

From the last lemma directly follows that (\ref{epsilon1}) is equal to
\[
(-1)^{m}\dfrac{s(2m+2s+1)}{(m+s)(m+s+1)}{m+2s \choose 2s}{2s \choose s},
\]
and (\ref{epsilon}) becomes
\begin{equation}\label{epsilonf}
\epsilon_{n,k}^{(s)}=\left\{\begin{array}{ll}
\pi(-1)^{m}\displaystyle\dfrac{s(2m+2s+1)}{(m+s)(m+s+1)}{m+2s \choose 2s}{2s \choose s}, & k=2nm, \ m\in\mathbb{N}_0, \\[0.in]
 0, & \mbox{\rm{otherwise}}. \\
\end{array}\right.
\end{equation}
Using the obtained
results, we get
\begin{eqnarray*}
|R_{n,s}(f)|&=&\left|\sum_{k=0}^{+\infty}\alpha_{(2s+1)n+k}\epsilon_{n,k} \right|=\left|\sum_{k=0}^{+\infty}\alpha_{(2s+1)n+2jn}\epsilon_{n,2jn} \right|\\
&\leq&\dfrac{2\pi}{\rho^{(2s+1)n}}\left(\max_{z\in
\mathcal{E}_{\rho}}|f(z)|\right)\sum_{m=0}^{+\infty}{\frac{s(2m+2s+1)}{(m+s)(m+s+1)}{m+2s \choose 2s}{2s \choose s}}{\rho^{-2mn}}\\
&=&2\pi\rho^{-n}\left(\max_{z\in \mathcal{E}_{\rho}}|f(z)|\right)F(x),
\end{eqnarray*}
where $x=\rho^{-n}$ (hence, $x\in(0,1)$) and
\[
F(x)=s{2s \choose s}\sum_{m=0}^{+\infty}{m+2s \choose
2s}\frac{(2m+2s+1)x^{m+s}}{(m+s)(m+s+1)}\,.
\]

\begin{lemma} If $x\in(0,1)$, then the sum of the series F(x) is equal to
\[
\dfrac{\sum_{k=0}^s(-1)^{k}{2s+1 \choose s+k+1}x^{s+k}}{(1-x)^{2s}}.
\]
\end{lemma}

\noindent{\bf Proof.}
We have that
\begin{eqnarray*}
F(x)&=&{\sum_{k=0}^s(-1)^{k}{2s+1 \choose s+k+1}x^{s+k}}\cdot{(1-x)}^{2s}\\
&=&{\sum_{k=0}^s(-1)^{k}{2s+1 \choose s+k+1}x^{s+k}}\sum_{j=0}^{+\infty}{j+2s-1 \choose 2s-1}x^{j}\\
&=&x^s{\sum_{k=0}^s(-1)^{k}{2s+1 \choose s+k+1}x^{k}}\sum_{j=0}^{+\infty}{j+2s-1 \choose 2s-1}x^{j}.
\end{eqnarray*}
With the aim of showing that the last is equal to
\[
x^s \cdot s{2s \choose s}\sum_{m=0}^{+\infty}{m+2s \choose
2s}\frac{(2m+2s+1)x^{m}}{(m+s)(m+s+1)},
\]
we actually have to show
\[
{{\sum_{k=0}^s}(-1)^{k}{2s+1 \choose s+k+1}}{j+2s-1 \choose 2s-1}=s{2s \choose s}\sum_{m=0}^{+\infty}{m+2s \choose
2s}\frac{(2m+2s+1)x^{m}}{(m+s)(m+s+1)}
\]
under the condition $k+j=m$, i.\,e.
\[
{{\sum_{j=m-s}^m}(-1)^{m-j}{2s+1 \choose s-(m-j)}}{j+2s-1 \choose j}=s{2s \choose s}\sum_{m=0}^{+\infty}{m+2s \choose
2s}\frac{(2m+2s+1)x^{m}}{(m+s)(m+s+1)},
\]
which is again the relation between (\ref{epsilon}) and (\ref{epsilonf}).\hfill\qed

Finally, we can formulate the main result in this part.
\begin{theorem}
If the function $f$ is analytic in the interior of the region $\mathcal{D}$ bounded by the $\mathcal{E}_{\rho}$ and continuous on $\overline{\mathcal{D}}$, the following error bound holds
\begin{equation}\label{ocenan}
\left|R_{n,s}(f)\right|\leq 2\pi\left(\max_{z\in
\mathcal{E}_{\rho}}|f(z)|\right)\dfrac{\sum_{k=0}^s(-1)^{k}{2s+1 \choose s-k}\rho^{2n(s-k)}}{\rho^n(\rho^{2n}-1)^{2s}}\,.
\end{equation}
\end{theorem}

\section{$L^1$-error bounds}
According to (\ref{eocena}) we study now the quantity
\[
L_{n,s}(\mathcal{E}_{\rho})=\dfrac{1}{2\pi}\oint_{\mathcal{E}_{\rho}}\left|K_{n,s}(z)\right||dz|,
\]
where $|K_{n,s}(z)|$ can be obtained from (\ref{kvadratmodula}). Since
$z=(u+u^{-1})/{2}$, $u=\rho e^{i\theta}$, and
$|dz|=(1/\sqrt{2})\cdot\sqrt{a_2-\cos{2\theta}}\,d\theta$ (see
\cite{Hunter}), the quantity $L_{n,s}(\mathcal{E}_{\rho})$ reduces to
\begin{equation}\label{izraz}
\begin{array}{rl}
L_{n,s}(\mathcal{E}_{\rho})=&\displaystyle\dfrac{1}{2\pi\sqrt{2}}\int_{0}^{2\pi}|K_{n,s}(z)|\sqrt{a_2-\cos{2\theta}}\,d\theta\\
=&\displaystyle\dfrac{1}{2\sqrt{2}}\int_{0}^{2\pi}\dfrac{\sqrt{a}}{\rho^{(2s+1)n}2^{s-1/2}c^s}\,d\theta=
\dfrac{1}{2^s\rho^{(2s+1)n}}\int_0^{\pi}\dfrac{\sqrt{a}}{c^s}\,d\theta.\\
\end{array}
\end{equation}
Applying Cauchy inequality to the last expression, we obtain
\begin{equation}\label{Locena}
L_{n,s}(\mathcal{E}_{\rho})\leq\dfrac{\sqrt{\pi}}{2^s \rho^{(2s+1)n}}
\sqrt{\int_0^{\pi}\dfrac{a}{c^{2s}}\,d\theta},
\end{equation}
where $a,c$ are given in (\ref{kvadratmodula}).

We have
\begin{eqnarray*}
a&=&\displaystyle\sum_{k=0}^s {2s+1\choose k}^2\rho^{2nk}+2\sum_{i<j, \, 0\leq i,j\leq s}{2s+1\choose i} {2s+1\choose j}\rho^{{2n}(i+j)}\cos{2n(j-i)\theta}\\
&=&\displaystyle\sum_{k=0}^s {2s+1\choose k}^2\rho^{2nk}+2\sum_{l=1}^{s}\cos{2nl\theta}\sum_{i=0}^{s-l}{2s+1\choose i} {2s+1\choose i+l}\rho^{{2n}(2i+l)},
\end{eqnarray*}
and
\begin{equation}\label{integral}
\int_0^{\pi}\dfrac{a}{c^{2s}}\,d\theta=\sum_{k=0}^s {2s+1\choose k}^2\rho^{4nk}I_0+2\sum_{l=1}^{s}\rho^{2nl}\sum_{i=0}^{s-l}{2s+1\choose i} {2s+1\choose i+l}\rho^{4ni}I_l,
\end{equation}
where
\begin{equation}\label{I0}
I_0=\int_0^{\pi}\dfrac{d\theta}{(a_{2n}+\cos{2n\theta})^{2s}}=(2\rho^{2n})^{2s}\dfrac{\pi\sum_{m=0}^{2s-1}{2s-1 \choose m}{4s-m-2 \choose 2s-1}\left(\rho^{4n}-1\right)^m}{\left(\rho^{4n}-1\right)^{4s-1}}
\end{equation}
and
\begin{equation}\label{I1}
I_l=\int_0^{\pi}\dfrac{\cos{2nl\theta}\,d\theta}{(a_{2n}+\cos{2n\theta})^{2s}}=(-1)^l(2\rho^{2n})^{2s}\dfrac{\pi\sum_{m=0}^{2s-1}{2s+l-1 \choose m}{4s-m-2 \choose 2s-1}\left(\rho^{4n}-1\right)^m}{\rho^{2nl}\left(\rho^{4n}-1\right)^{4s-1}}.
\end{equation}
We have used \cite[Eq. 3.616.7]{Ryzhik}.

\section{Numerical example}\label{sect6}

We consider the calculation of the integral
\[
I_\omega(f)=\int_{-1}^1 f(t)T_n(t)/\sqrt{1-t^2}\,dt
\]
by using the quadrature formula (\ref{glavnaqf}), where
the function
\[
f(z)=f_0(z)=e^{\omega z^2}\quad(\omega>0)
\]
 is entire. We tested the derived bounds for some values of $n$, $s$ and $\omega>0$.
 Since $f_0$ is an entire function, the different estimations hold for $\mathcal{E}_\rho$, with $\rho\in(1,\infty)$. It is easy to see that
\[
\max_{z\in\mathcal{E}_\rho}\left|e^{\omega z^2}\right|=e^{\omega a_1^2},\
a_1={\textstyle\frac 12}(\rho+\rho^{-1}) .
\]
 \begin{table}[htb]
	\begin{center}{\footnotesize
			\begin{tabular}{llllll}\hline
				$n,s,\omega$ & $r_1(f_0)$ & $r_2(f_0)$
				& $r_3(f_0)$ & ${\rm Error}$ & $I_\omega$ \\
				\hline
				$8, 1, 1$ & $5.22(-14) $ & $3.40(-14) $ & $1.70(-14)$ & $ 1.94(-15)$ &
				$8.53...(-4)$ \\
				$8, 2, 1$ & $7.36(-28) $ & $4.38(-28) $ & $2.19(-28) $ & $ 1.94(-29)$ &
				$8.53...(-4)$ \\
				$8, 3, 1$ & $3.95(-43) $ & $2.20(-43) $ & $1.10(-43)$ & $ 8.27(-45)$ &
				$8.53...(-4)$ \\
				$8, 1, 5$ & $9.05(-5) $ & $6.95(-5) $ & $3.47(-5) $ & $ 3.93(-7)$ &
				$5.28...(+0)$ \\
				$8, 2, 5$ & $4.43(-13) $ & $3.32(-13) $ & $1.66(-13) $ & $ 1.47(-14)$ &
				$5.28...(+0)$ \\
				$8, 3, 5$ & $8.83(-23) $ & $6.41(-23) $ & $3.12(-23)$ & $ 2.40(-24)$ &
				$5.28...(+0)$ \\
				$8, 1, 10$ & $7.02(+0) $ & $5.10(+0) $ & $2.55(+0) $ & $2.79(-1)$ &
				$2.38...(+3)$ \\
				$8, 2, 10$ & $7.07(-6) $ & $5.41(-6) $ & $2.71(-6) $ & $2.34(-7)$ &
				$2.38...(+3)$ \\
				$8, 3, 10$ & $3.26(-13) $ & $2.49(-13) $ & $1.24(-13) $ & $9.21(-15)$ &
				$2.38...(+3)$ \\
				$8, 1, 20$ & $2.14(+7) $ & $1.30(+7) $ & $1.51(+6) $ & $6.46(+5)$ &
				$8.48...(+7)$ \\
				$8, 2, 20$ & $2.94(+3) $ & $2.09(+3) $ & $1.05(+3) $ & $8.66(+1)$ &
				$8.48...(+7)$ \\
				$8, 3, 20$ & $2.54(-2) $ & $1.90(-2) $ & $ 9.50(-3)$ & $6.92(-4)$ &
				$8.48...(+7)$ \\
				\hline
				$12, 1, 1$ & $1.25(-24) $ & $7.57(-25) $ & $3.79(-25) $ & $3.54(-26)$ &
				$1.77...(-6)$ \\
				$12, 2, 1$ & $8.51(-48) $ & $4.68(-48) $ & $2.34(-48)$ & $ 1.70(-49)$ &
				$1.77...(-6)$ \\
				$12, 3, 1$ & $4.24(-73) $ & $2.18(-73) $ & $1.09(-73)$ & $6.68(-75)$ &
				$1.77...(-6)$ \\
				$12, 1, 5$ & $3.07(-11) $ & $2.32(-11) $ & $1.16(-11) $ & $1.08(-12)$ &
				$2.52...(-1)$ \\
				$12, 2, 5$ & $4.70(-26) $ & $3.38(-26) $ & $1.69(-26) $ & $ 1.23(-27)$ &
				$2.52...(-1)$ \\
				$12, 3, 5$ & $5.51(-43) $ & $3.79(-43) $ & $1.90(-43)$ & $1.16(-44)$ &
				$2.52...(-1)$ \\
				$12, 1, 10$ & $1.26(-4) $ & $9.59(-5) $ & $4.79(-5) $ & $4.41(-6)$ &
				$3.69...(+2)$ \\
				$12, 2, 10$ & $6.80(-16) $ & $5.18(-16) $ & $2.59(-16) $ & $1.86(-17)$ &
				$3.69...(+2)$ \\
				$12, 3, 10$ & $3.04(-29) $ & $2.27(-29) $ & $1.84(-29) $ & $6.94(-31)$ &
				$3.69...(+2)$ \\
				$12, 1, 20$ & $1.48(+4) $ & $1.02(+4) $ & $5.12(+3) $ & $4.43(+2)$ &
				$3.10...(+7)$ \\
				$12, 2, 20$ & $2.03(-4) $ & $1.53(-4) $ & $7.64(-5) $ & $5.39(-6)$ &
				$3.10...(+7)$ \\
				$12, 3, 20$ & $3.02(-14) $ & $2.31(-14) $ & $1.16(-14) $ & $6.98(-16)$ &
				$3.10...(+7)$ \\
				\hline
				$16, 1, 1$ & $3.83(-36) $ & $2.20(-36) $ & $1.10(-36) $ & $8.92(-38)$ &
				$1.97...(-9)$ \\
				$16, 2, 1$ & $3.23(-69) $ & $1.67(-69) $ & $8.37(-70)$ & $ 5.26(-71)$ &
				$1.97...(-9)$ \\
				$16, 3, 1$ & $3.83(-105) $ & $1.85(-105)$ & $9.24(-106)$ & $4.92(-107)$ &
				$1.97...(-9)$ \\
				$16, 1, 5$ & $1.40(-18) $ & $1.03(-18) $ & $5.16(-19) $ & $4.17(-20)$ &
				$6.72...(-3)$ \\
				$16, 2, 5$ & $1.69(-40) $ & $1.17(-40) $ & $5.85(-41)$ & $ 3.67(-42)$ &
				$6.72...(-3)$ \\
				$16, 3, 5$ & $2.96(-65) $ & $1.94(-65) $ & $9.73(-66)$ & $5.18(-67)$ &
				$6.72...(-3)$ \\
				$16, 1, 10$ & $3.34(-10) $ & $2.56(-10) $ & $1.28(-10) $ & $1.03(-11)$ &
				$3.45...(+1)$ \\
				$16, 2, 10$ & $2.34(-27) $ & $1.76(-27) $ & $8.78(-28) $ & $5.50(-29)$ &
				$3.45...(+1)$ \\
				$16, 3, 10$ & $2.58(-47) $ & $1.87(-47) $ & $9.35(-48)$ & $4.93(-49)$ &
				$3.45...(+1)$ \\
				$16, 1, 20$ & $1.90(+0) $ & $1.38(+0) $ & $6.9(-1) $ & $5.32(-2)$ &
				$8.03...(+6)$ \\
				$16, 2, 20$ & $6.10(-13)$ & $4.67(-13) $ & $2.34(-13) $ & $1.41(-14)$ &
				$8.03...(+6)$ \\
				$16, 3, 20$ & $3.69(-28) $ & $2.82(-28) $ & $1.41(-28)$ & $7.32(-30)$ &
				$8.03...(+6)$ \\
				\hline
		\end{tabular}}
		\vskip 0.2in
		\caption{The values of the derived bounds $r_1(f_0),r_2(f_0),r_3(f_0)$, the actual (sharp) errors, and the values of the integrals $I_\omega$, for some values of $n$,  $s$, $\omega$.}\label{figura111}
	\end{center}
\end{table}

The length of the ellipse $\mathcal{E}_\rho$ can be estimated by
(cf. \cite[Eq. (2.2)]{scherer})
\begin{equation}\label{duzinaelipse}
\ell(\mathcal{E}_\rho)\le2\pi a_1\left(1-\frac 14
a_1^{-2}-\frac{3}{64}a_1^{-4}-\frac{5}{256}a_1^{-6}\right).
\end{equation}
The corresponding bounds ($|R_{n,s}(f)|\le r_i(f), i=1,2,3$) have the form
\[ \
r_1(f)=\inf_{\rho_0<\rho<+\infty} B_1, \quad r_2(f)=\inf_{1<\rho<+\infty}B_2,\quad r_3(f)=\inf_{1<\rho<+\infty} B_3,
\]
($\rho_0$ is defined in Theorem \ref{teorema1})
where from (\ref{ocena}) and (\ref{kvadratmodula})
\[
B_1=\dfrac{\pi a_1}{2^{s-\frac{1}{2}}\rho^{(2s+1)n}}\cdot\dfrac{{\sum_{k=0}^{s}{2s+1 \choose k}\rho^{2nk}}}{\sqrt{a_2-1}(a_{2n}+1)^{s}}\left(1-\frac 14 a_1^{-2}-\frac{3}{64}a_1^{-4}-\frac{5}{256}a_1^{-6}\right)e^{\omega a_1^2},
\]
from (\ref{ocenan})
\[
B_2=2\pi\dfrac{\sum_{k=0}^s(-1)^{k}{2s+1 \choose s-k}\rho^{2n(s-k)}}{\rho^n(\rho^{2n}-1)^{2s}} \,e^{\omega a_1^2},
\]
and from (\ref{eocena}), (\ref{Locena}) and (\ref{integral})
\begin{eqnarray*}
B_3&=&\dfrac{\sqrt{\pi}}{2^s \rho^{(2s+1)n}}
\sqrt{\sum_{k=0}^s {2s+1\choose k}^2\rho^{4nk}I_0+2\sum_{l=1}^{s}\rho^{2nl}\sum_{i=0}^{s-l}{2s+1\choose i} {2s+1\choose i+l}\rho^{4ni}I_l}\\
&&\displaystyle\,\times\, e^{\omega a_1^2}.
\end{eqnarray*}
$I_0$, $I_l$ are given by (\ref{I0}), (\ref{I1}) respectively, and $a_j$ by (\ref{aiovi}). The corresponding results are displayed in Table \ref{figura111}. In Table \ref{figura111} are also displayed the actual (sharp) errors ``Error'' and the values of the integrals $I_\omega$.

On the basis of displayed results in  Table \ref{figura111} we conclude the all three kind of considered error bounds are of the same range, they are also very close to the actual error. In order to find a quadrature sum $Q_{n,s}(f)=\sum_{\nu=1}^n\sum_{i=0}^{2s-1}A_{i,\nu}f^{(i)}(\xi_\nu)$ in (\ref{glavnaqf}) one has to calculate $2sn$ values $f^{(i)}(\xi_\nu)$. An error bound of $Q_{n,s}(f)$ in  Table \ref{figura111} is of the form $C\cdot 10^{-l}$ $(1\le C<10$). It is clear from Table \ref{figura111} that if we fix $\omega$ (the integrand) and $n$ (the number of nodes), then the error bounds of the same kind decrease. So, if we instead of $Q_{n,s}(f)$ calculate $Q_{n,s+1}(f)$, the amount of computations of $f^{(i)}(\xi_\nu)$ increases in $2n$, and the corresponding error bound decreases in dependance on the integrand.

\begin{acknowledgements} The authors are indebted to the unknown referees for the valuable comments that have improved the first version of the paper.
Research supported in part by the Serbian Ministry of Education, Science and Technological Development (Research Project:
``Methods of numerical and nonlinear analysis with applications'' (\# 174002)).
\end{acknowledgements}

%\vfill\break

% BibTeX users please use
%\bibliographystyle{spmpsci}

\begin{thebibliography}{33}

\bibitem{boj} Bojanov, B.: On a quadrature formula of Micchelli and
Rivlin. J. Comput. Appl. Math. {\bf 70}, 349--356 (1996)
\bibitem{bojanov} Bojanov, B., Petrova, G.: Quadrature formulae for
Fourier coefficients. J. Comput. Appl. Math. {\bf 231},  378--391 (2009)
\bibitem{rdevore} DeVore, R.: A property of Chebyshev polynomials.
J. Approx Theory {\bf 12},  418--419 (1974)
\bibitem{Eliot} Elliot, D.: The evaluation and estimation of the coefficients in the Chebyshev series
expansion of a functions. Math. Comp. {\bf 18},  82--90 (1964)
\bibitem{gauvar} Gautschi, W., Varga, R.S.: Error bounds for Gaussian
quadrature of analytic functions. SIAM J. Numer. Anal. {\bf 20},
1170--1186  (1983)
\bibitem{gon} Gon\v carov, V.L.: Theory of of interpolation and
approximation of functions. GITTL, Moscow, 1954 (in Russian)
\bibitem{Ryzhik} Gradshteyn, I.S., Ryzhik, I.M.: Tables of integrals, series and products. 6th edn
(Jeffrey, A., Zwillinger, D., eds), Academic Press, San Diego,
2000
\bibitem{Hunter} Hunter, D.B.: Some error expansions for Gaussian quadrature. BIT Numer. Math. {\bf 35}, 64--82  (1995)
\bibitem{geno} Hunter, D.B., Nikolov, G.: On the error
term of symmetric Gauss-Lobatto quadrature formulae for analytic
functions. Math. Comp. {\bf 69}, 269--282  (2000)
\bibitem{mr} Micchelli, C.A.,  Rivlin, T.J.:  Tur\' an formulae and
highest precision quadrature rules for Chebyshev coefficients. IBM
J. Res. Develop. {\bf 16},  372--379 (1972)
\bibitem{mr10} Micchelli, C.A., Rivlin, T.J.: Some new characterizations of the Chebyshev polynomials. J. Approx. Theory {\bf 12},  420--424 (1974)
\bibitem{Baza} Milovanovi\'{c}, G.V., Spalevi\'{c}, M.M.:
An error expansion for some Gauss-Tur\'{a}n quadratures and
$L^1$-estimates of the remainder term. BIT Numer. Math.
{\bf 45},  117--136 (2005)
\bibitem{MSmc2013} Milovanovi\'{c}, G.V., Spalevi\'{c}, M.M.:
Kronrod extensions with multiple nodes of quadrature formulas for Fourier coefficients. Math. Comp.
{\bf 83},  1207--1231 (2014)
\bibitem{GVM_RO_MMS} Milovanovi\'{c}, G.V., Orive, R., Spalevi\'{c}, M.M.:
Quadrature with multiple nodes for Fourier-Chebyshev coefficient. IMA J. Numer. Anal., to appear. DOI: 10.1093/imanum/drx067
\bibitem{filomat} Milovanovi\'{c}, G.V., Pej\v{c}ev, A.V., Spalevi\'{c}, M.M.: A note on an error bound of Gauss-Tur\' an
quadrature with the Chebyshev weight. FILOMAT {\bf 27}, 1037-1042 (2013)
\bibitem{APejcevMMS} Pej\v cev, A.V., Spalevi\' c, M.M.: Error bounds of Micchelli-Rivlin quadrature formula for analytic functions. J. Approx. Theory {\bf 169},  23--34 (2013)
\bibitem{NM2016PejSpa} Pej\v cev, A.V., Spalevi\' c, M.M., The error bounds of Gauss-Radau quadrature formulae with Bernstein-Szeg\H{o} weight functions. Numer. Math. {\bf 133},  177--201 (2016)
\bibitem{scherer} Scherer, R.,  Schira, T.: Estimating quadrature errors
for analytic functions using kernel representations and
biorthogonal systems. Numer. Math. {\bf 84},  497--518 (2000)
\bibitem{sch} Schira, T.:  The remainder term for analytic functions of symmetric Gaussian quadratures.
 Math. Comp. {\bf 66}, 297--310  (1997)

\end{thebibliography}


\end{document}